\documentclass[11pt]{article}
\usepackage{amssymb,amsfonts,amsmath,latexsym,epsf,tikz,url}
\newtheorem{theorem}{Theorem}[section]

\newtheorem{observation}[theorem]{Observation}

\newcommand{\proof}{\noindent{\bf Proof.\ }}
\newcommand{\qed}{\hfill $\square$\medskip}
\allowdisplaybreaks

\begin{document}

\title{More on the  sixth coefficient of the matching
	polynomial in regular graphs}

\author{Neda Soltani \and Saeid Alikhani$^{}$\footnote{Corresponding author} }

\date{\today}

\maketitle

\begin{center}

Department of Mathematics, Yazd University, 89195-741, Yazd, Iran\\
{\tt  neda\_soltani@ymail.com, alikhani@yazd.ac.ir }\\

\end{center}


\begin{abstract}
A matching set $M$ in a graph $G$ is a collection of edges of $G$ such that no two edges from $M$ share a vertex.
In this paper we consider some parameters related to the matching of regular  graphs.
We find the sixth coefficient of the matching polynomial  of 
regular graphs. As a consequence,  every cubic graph of order $10$  is matching unique.
\end{abstract}


\medskip
\noindent{\bf AMS Subj.\ Class.:} 05C30; 05C70.

\section{Introduction}

A set of independent edges of a graph $G$ is called a matching of $G$. A matching $M$ is maximum if there is no matching 
in $G$ with more edges than $M$. The cardinality of any maximum matching in $G$ 
is called the matching number of $G$ and is denoted by $\alpha^{\prime}(G)$. If each vertex of $G$ is incident with an edge of $M$, the matching $M$ is called perfect. Only graphs of even order $n$ can have a perfect matching and the size of such a matching is $\frac{n}{2}$. A matching $M$ in $G$ is maximal if no other matching in $G$ contains it as a proper subset.  
The cardinality of any smallest maximal matching in $G$ is the saturation number of $G$ and is denoted by $s(G)$ (the same term, saturation number, is also used in the literature with a different meaning; we refer the reader to \cite{Fad} for more information).

	If a graph $G$ has a maximum matching of size $k$, then any maximal matching has at least size $\frac{k}{2}$ (\cite{Bie}).
This  implies that $s(G)\geq \frac{\alpha^{\prime}(G)}{2}$. 
We recall that a set of vertices $I$ is independent if no two vertices from $I$ are adjacent. Clearly, the set of vertices that is not covered by a maximal matching is independent. This observation provides an obvious lower bound on saturation number of the graph $G$, i.e. $s(G)\geq \frac{(n-|I|)}{2}$ where $G$ is graph of order $n$ (\cite{Ves}).

 An $r$-matching in a graph $G$ is a  set of $r$ pairwise non-incident edges. The number of $r$-matchings in $G$ is denoted by $\rho(G,r)$. The
 matching polynomial of $G$ is defined by
$$ \mu(G,x)=\sum_{r=0}^{\lfloor\frac{n}{2}\rfloor} (-1)^r \rho(G,r)x^{n-2r},$$
 where $n$ is the order of $G$ and $p(G,0)$ is considered to be $1$ (see \cite{Gutman1,Gutman2,Gutman3,Gutman4}). Two graphs are co-matching if their matching polynomials are equal. A graph that
 is characterized by its matching polynomial is said to be matching unique. 
Computation of matching polynomial is equivalent to computation of  the number of $k$-matchings in the
graph, for all $k$.  The Petersen graph is one of the famous graph which is a symmetric
non-planar cubic graph of order $10$.  Behmaram established a formula for the number of $4$-matchings in triangular-free graph with respect to the
number of vertices, edges, degrees and $4$-cycles. Also he has proved that  the Petersen
graph is uniquely determined by its matching polynomial \cite{Beh}. The Petersen graph is one of the cubic graphs of order $10$ and this is natural to study the matching polynomial of all of cubic graphs of order $10$ which some of them are not triangular free. 
Vesalian and Asgari in \cite{Match} has obtained a formula for the number of $5$-matchings in triangular-free and $4$-cycle-free
graph based on the number of vertices, edges, the degrees of vertices and the number of $5$-cycles.


\medskip

In the next section, after investigation of the saturation number and the matching number of cubic graphs of order $10$, we  establish a formula for the number of $5$-matchings in regular graphs. Using our result, we present the matching polynomial of cubic graphs of order $10$ and conclude that every cubic graph of order  ten  is matching unique.
\section{Main results} 

First  we consider  the saturation number of cubic graphs of order $10$ as a special kind of regular graphs. We recall that a cubic graph is a $3$-regular graph. There are exactly $21$ cubic graphs of order $10$ given in Figure 1 (see \cite{Kho}). Note that the graph $G_{17}$ is the Petersen graph.
  The following observation presents the saturation number of cubic graphs with $10$ vertices. Note that $s(G_{21})=5$.

\begin{figure}[h]
	\vglue-1cm
	\hglue2.5cm
	\includegraphics[width=11cm,height=4.9cm]{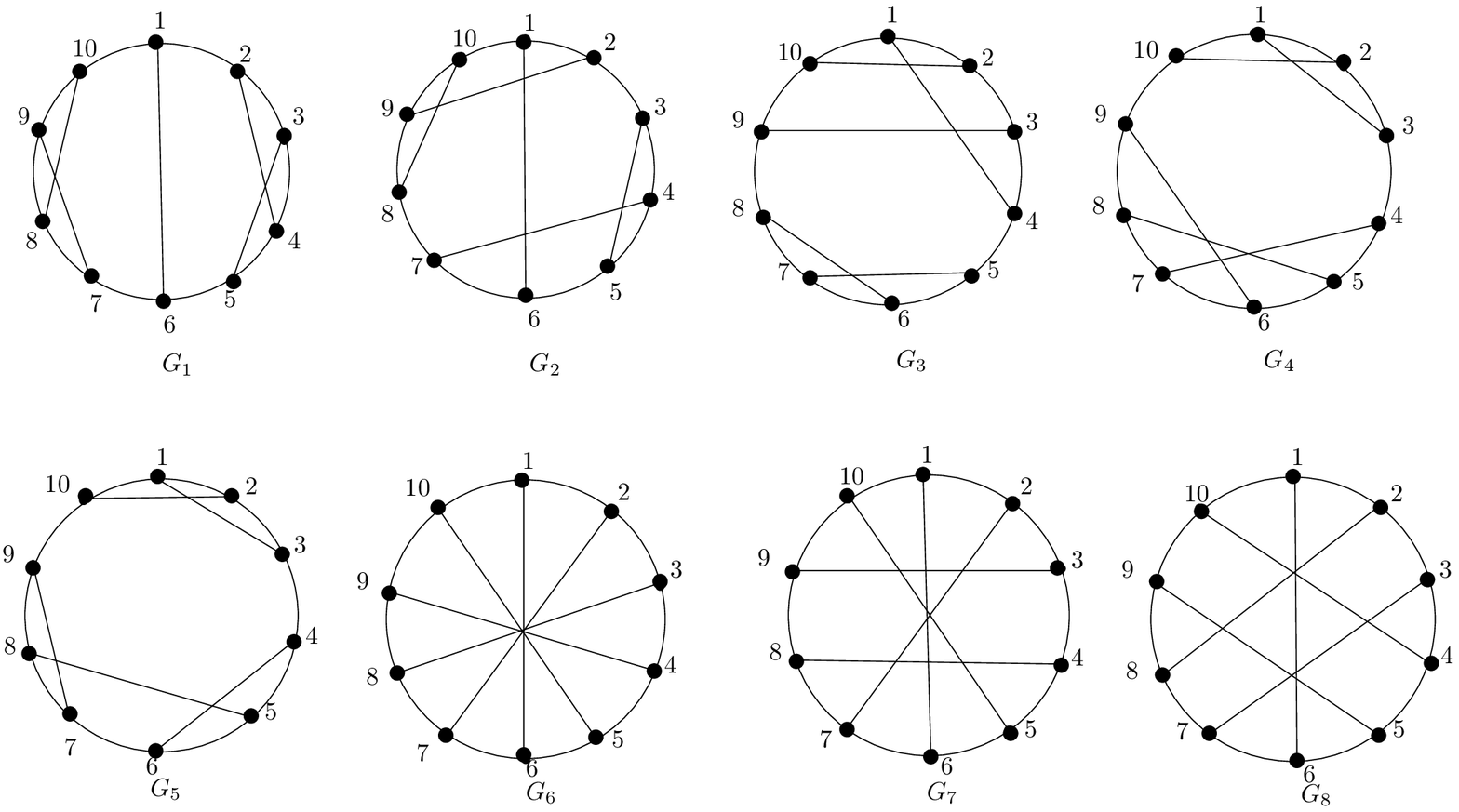}
	\vglue5pt
	\hglue2.5cm
	\includegraphics[width=11cm,height=4.9cm]{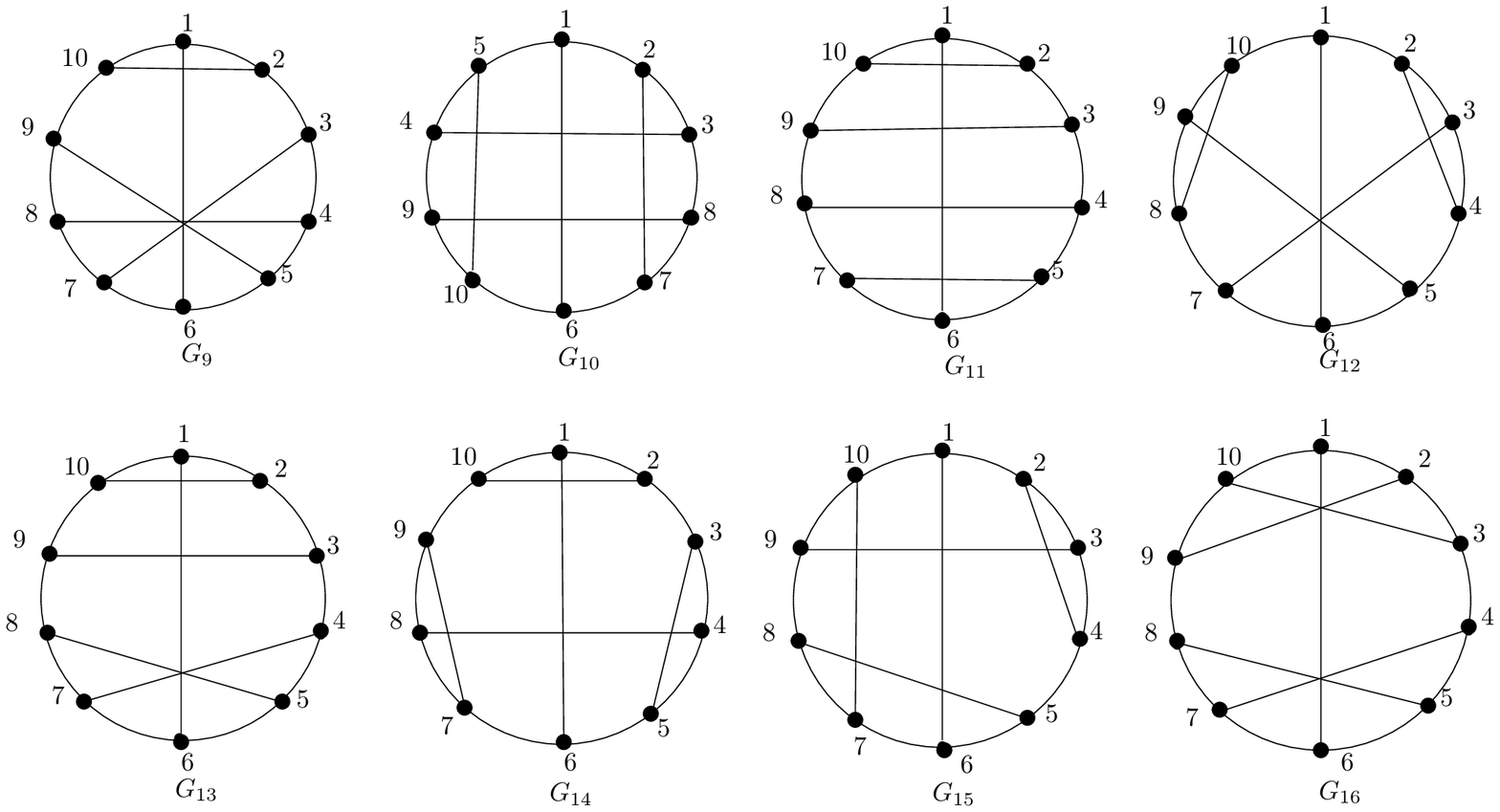}
	\hglue2.5cm
	\includegraphics[width=10.7cm,height=4.9cm]{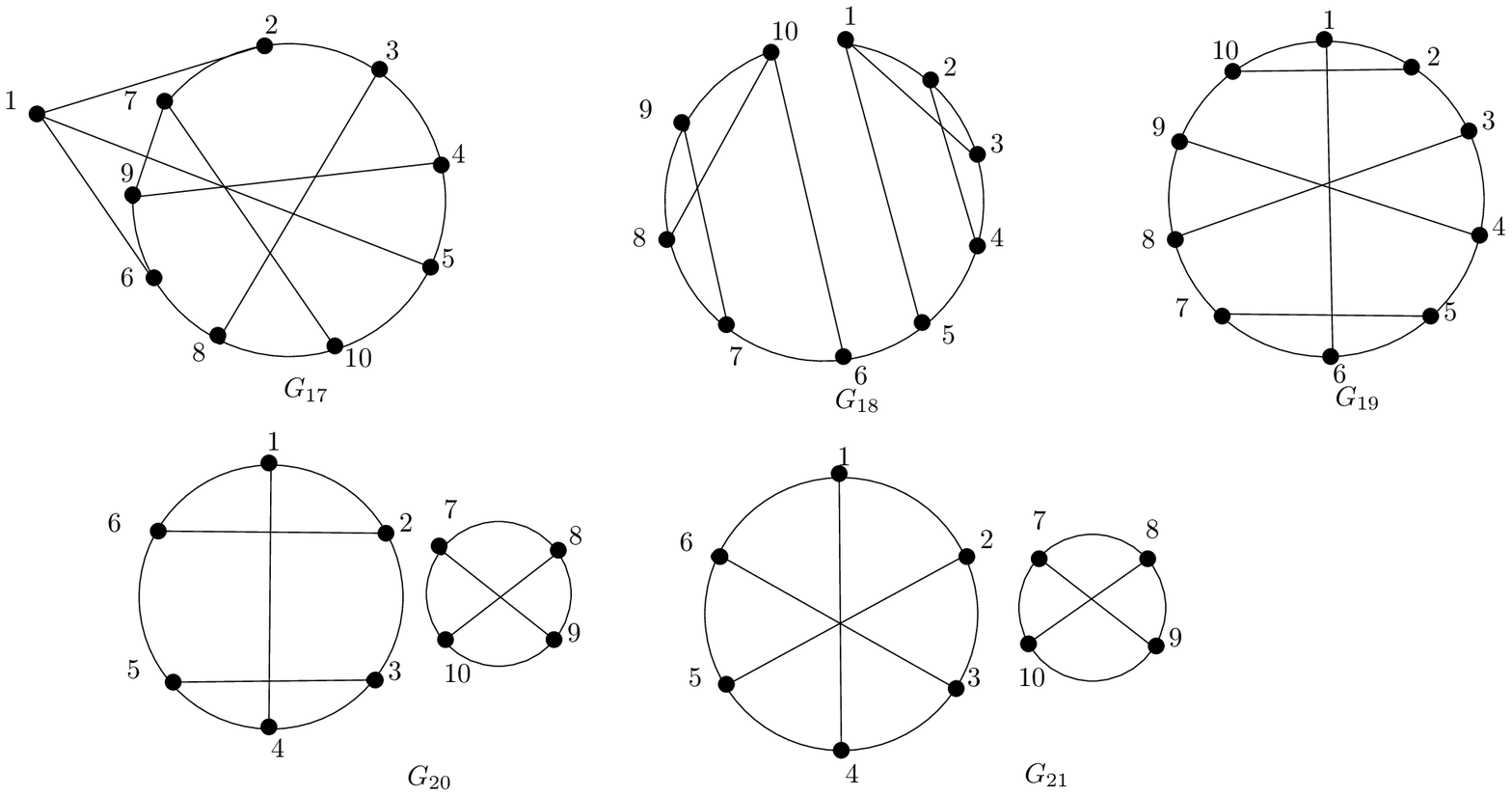}
	\vglue-10pt \caption{\label{figure1} \small {Cubic graphs of order $10$}.}
\end{figure}

\begin{observation} 
	\begin{enumerate}
		\item [(i)]
		$s(G_i)=3$ for $i=1,2,5,8,12,16,17,18.$
		\item[(ii)]
		$s(G_i)=4$ for $i=3,4,6,7,9,10,11,13,14,15,19,20.$
		\end{enumerate}
\end{observation} 

\medskip 

Since each $2$-connected $3$-regular graph has a perfect matching  and a maximum matching in a disconnected graph consists of the union of maximum matchings in each of its components, so we can conclude that every cubic graph 
of order $10$ has perfect matching. Thus for each $1\leq i\leq 21$, $\alpha^{\prime}(G_i)=5$.  

\medskip
The domination polynomial  and the edge cover
polynomial of cubic graphs of order $10$ has studied in \cite{turk} and \cite{jas}, respectively. It has proved that  all cubic
graphs of order $10$  are determined
uniquely by their edge cover polynomials, but this is not true for their domination polynomials. 
Here we consider  the matching polynomial of regular graphs. Suppose that $G$ is a regular graph of degree $r$ with $n$ vertices and $\Gamma_{m, k, i}$ is an isomorphism class containing $i$ kind of spanning subgraphs of $G$, 
with $m$ edges from $G$ and $k$ vertices of degree $1$. Let $G_{m, k, i}$ denote a representative chosen from 
$\Gamma_{m, k, i}$ and  $g_{m,k,i}$ is the cardinality of the set $\Gamma_{m, k, i}$ (\cite{bee2}). Then the 
matching polynomial of graph $G$ can be expressed as follows
$$\mu(G, x)=\sum_{m=0}^{\left\lfloor\frac{n}{2}\right\rfloor}(-1)^mg_{m, 2m, 1}x^{n-2m}.$$
Note that $g_{m, 2m, 1}$ denotes the number of $m$-matchings of the graph $G$, for each $m$ and that $g_{0, 0, 1}=1$. The following theorem gives an approach  for acquiring the coefficients of the matching polynomial of regular graphs.

\begin{theorem}{\rm\cite{Bee}}\label{Ce}
Suppose that $G$ is a regular graph of degree $r$ on $n$ vertices. Then for each combination of $m$, $k$ and $i^*$ there exist constants $a_{j, i}=a_{j, i}(n, r)$, that depend on $G$ only through $n$ and $r$, and constants $a_i$ that are independent of $G$, so that $$g_{m, k, i^*}=\sum_{j=0}^{m-1}\sum_{i}a_{j, i}(n, r)g_{j, 0, i}+\sum_{i}a_ig_{m, 0, i}.$$
\end{theorem}

By applying the proof of Theorem \ref{Ce}, Beezer obtained the following equations for $g_{1, 2, 1}, g_{2, 2, 1}, ..., 
g_{4, 8, 1}$ in \cite{Bze}. 

\medskip
\medskip
\medskip
$\begin{array}{l}
\medskip
g_{1, 2, 1}=\frac{nr}{2},  \\ 
\medskip
g_{2, 2, 1}=\frac{n(r-1)r}{2},  \\ 
\medskip
g_{2, 4, 1}=\frac{nr}{8}(nr-4r+2),  \\ 
\medskip
g_{3, 2, 1}=\frac{n(r-1)^2r}{2}-3g_{3, 0, 1},  \\
\medskip
g_{3, 3, 1}=\frac{n(r-2)(r-1)r}{6},  \\
\medskip
g_{3, 4, 1}=\frac{n(r-1)r(nr-6r+4)}{4}+3g_{3, 0, 1},  \\
\medskip
g_{3, 6, 1}=\frac{nr}{48}(n^2r^2-12nr^2+40r^2+6nr-48r+16)-g_{3, 0, 1},  \\ 
\medskip
g_{4, 1, 1}=(3r-6)g_{3, 0, 1},  \\
\medskip
g_{4, 2, 1}=\big(\frac{nr}{2}-3r+3\big)g_{3, 0, 1},  \\
\medskip
g_{4, 2, 2}=\frac{n(r-1)^3r}{2}+(-6r+9)g_{3, 0, 1}-4g_{4, 0, 1}, \\
\medskip
g_{4, 3, 1}=\frac{n(r-2)(r-1)^2r}{2}+(-6r+12)g_{3, 0, 1},  \\
\medskip
g_{4, 4, 1}=\frac{n(r-1)^2r(nr-8r+6)}{4}+\big(-\frac{3nr}{2}+18r-21\big)g_{3, 0, 1}+4g_{4, 0, 1},  \\
\medskip
g_{4, 4, 2}=\frac{n(r-3)(r-2)(r-1)r}{24},  \\
\medskip
g_{4, 4, 3}=\frac{n(r-1)^2r(nr-9r+8)}{8}+(6r-9)g_{3, 0, 1}+2g_{4, 0, 1},  \\
\medskip
g_{4, 5, 1}=\frac{n(r-2)(r-1)r(nr-8r+6)}{12}+(3r-6)g_{3, 0, 1},  \\
\medskip
g_{4, 6, 1}=\frac{n(r-1)r(n^2r^2-16nr^2+72r^2+10nr-104r+40)}{16}+\big(\frac{3nr}{2}-21r+24\big)g_{3, 0, 1}
-4g_{4, 0, 1},  \\
\medskip
g_{4,8,1}=\frac{nr}{384}\big(n^3r^3-24n^2r^3+208nr^3-672r^3+12n^2r^2-240nr^2+1344r^2+76nr\\~~~~~~~-960r+240\big)+\big(-\frac{nr}{2}+6r-6\big)g_{3, 0, 1}+g_{4, 0, 1}.  
\end{array}$

\medskip
\medskip
\medskip
We follow  approach in \cite{Bze} to derive a formula for the sixth coefficient of the matching polynomial in regular graphs, i.e., $g_{5, 10, 1}$.  

\medskip
We start with possible subgraphs on five edges or less (see Figures 3 and 4) and determine the subgraphs with $5$ edges and a vertex of degree $1$. Then remove the edge incident to degree $1$ vertex and call other endpoint $w$. We identify vertices isomorphic to $w$ and then add back a single edge, attaching one end at vertices like $w$. Finally, we determine the types of subgraphs formed and the amount of overcounting. Note that the subgraphs with no degree $1$ vertices are free variables. For example, Let $H$ be a subgraph of regular graph $G$ with $5$ edges. 

\begin{figure}[ht]
\centerline{\includegraphics[width=5cm]{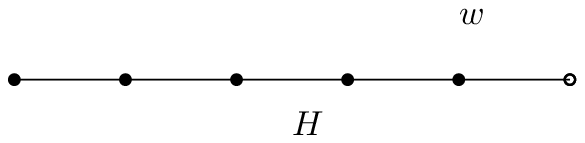}}
\end{figure}

Remove an edge incident to a vertex of degree $1$. Label other endpoint $w$. In the path on $4$ edges that remains, there is one other vertex like $w$. So there are 2 vertices like $w$ and $r-1$ ways to attach back an edge at $w$, considering all vertices as possibilities for the other end of the new edge (see Figure 2).

\begin{figure}[ht]
\centerline{\includegraphics[width=8.5cm]{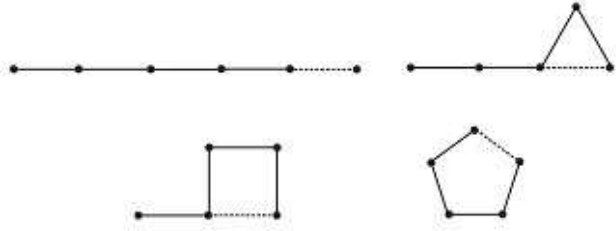}}
\caption{\label{5}\small {The subgraphs formed by adding back an edge at $w$}.}
\end{figure}

So we have 
$$2(r-1)g_{4, 2, 2}=2g_{5, 2, 2}+2g_{5, 1, 2}+2g_{5, 1, 1}+10g_{5, 0, 1}.$$

\begin{figure}[h]\label{sub}
	\hspace{2cm} 
	\includegraphics[width=10cm,height=8cm]{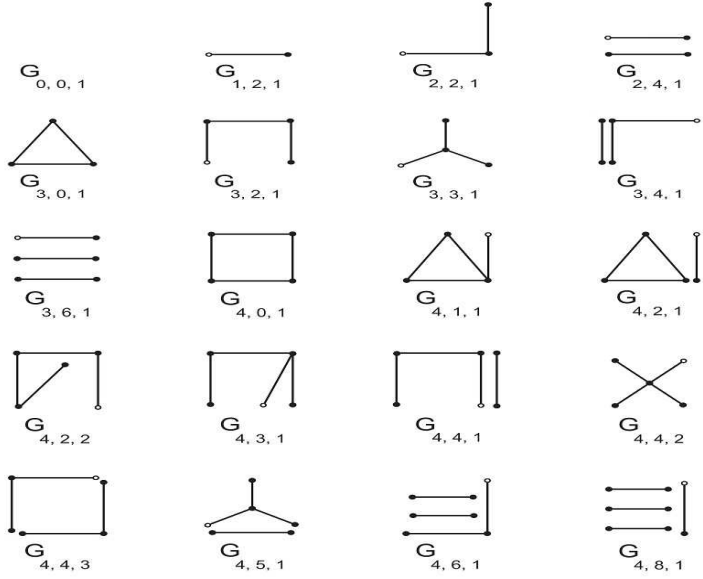}
      \caption{\label{5}\small {The possible subgraphs with $\leq 4$ edges. Vertex $w$ is adjacent to open-circle vertex}.}
\end{figure}

By using this process for every subgraph of regular graph $G$ with $5$ edges and a vertex of degree $1$, we have the following system of linear equations.

$$4(r-2)g_{4, 0, 1}=1g_{5, 1, 1}+2g_{5, 0, 2},$$
$$1(r-1)g_{4, 1, 1}=1g_{5, 1, 2}+4g_{5, 0, 2},$$ 
$$(n-4)rg_{4, 0, 1}=2g_{5, 2, 1}+1g_{5, 1, 1},$$
$$2(r-1)g_{4, 2, 2}=2g_{5, 2, 2}+2g_{5, 1, 2}+2g_{5, 1, 1}+10g_{5, 0, 1},$$
$$2(r-1)g_{4, 2, 1}=2g_{5, 2, 3}+1g_{5, 1, 2},$$
$$2(r-2)g_{4, 1, 1}=2g_{5, 2, 4}+4g_{5, 0, 2},$$
$$1(r-3)g_{4, 1, 1}=2g_{5, 2, 5},$$
$$2(r-2)g_{4, 2, 2}=2g_{5, 3, 1}+2g_{5, 1, 1}+2g_{5, 2, 4},$$
$$2(r-1)g_{4, 3, 1}=2g_{5, 3, 2}+2g_{5, 1, 2}+2g_{5, 2, 4}+2g_{5, 1, 1},$$
$$3(r-2)g_{4, 2, 1}=1g_{5, 3, 3}+1g_{5, 1, 2},$$
$$2(r-1)g_{4, 4, 1}=2g_{5, 4, 1}+2g_{5, 3, 3}+8g_{5, 2, 1}+2g_{5, 2, 2},$$
$$4(r-1)g_{4, 4, 3}=2g_{5, 4, 2}+6g_{5, 2, 3}+2g_{5, 2, 2}+1g_{5, 3, 1},$$
$$1(r-2)g_{4, 3, 1}=4g_{5, 4, 3}+2g_{5, 2, 4},$$
$$1(r-3)g_{4, 3, 1}=3g_{5, 4, 4}+2g_{5, 2, 5},$$
$$(n-5)rg_{4, 2, 1}=4g_{5, 4, 5}+2g_{5, 2, 3}+1g_{5, 3, 3},$$
$$3(r-1)g_{4, 5, 1}=1g_{5, 5, 1}+2g_{5, 3, 3}+1g_{5, 3, 1},$$
$$2(r-2)g_{4, 4, 3}=3g_{5, 5, 2}+1g_{5, 3, 1}+2g_{5, 4, 3},$$
$$1(r-4)g_{4, 4, 2}=5g_{5, 5, 3},$$
$$1(r-3)g_{4, 5, 1}=4g_{5, 6, 1}+1g_{5, 4, 4},$$
$$2(r-1)g_{4, 6, 1}=2g_{5, 6, 2}+6g_{5, 4, 5}+2g_{5, 4, 1},$$
$$4(r-1)g_{4, 6, 1}=4g_{5, 6, 3}+2g_{5, 4, 1}+1g_{5, 5, 1}+2g_{5, 4, 2},$$
$$1(r-2)g_{4, 6, 1}=3g_{5, 7, 1}+1g_{5, 5, 1},$$
$$8(r-1)g_{4, 8, 1}=2g_{5, 8, 1}+2g_{5, 6, 2},$$
$$(n-8)rg_{4, 8, 1}=10g_{5, 10, 1}+2g_{5, 8, 1}.$$

\medskip
\medskip
Now we can conclude that

\begingroup\makeatletter\def\f@size{5}\check@mathfonts
\def\maketag@@@#1{\hbox{\m@th\large\normalfont#1}}%
$\begin{array}{l}
\medskip

 g_{5, 1, 1}=(4r-8)g_{4, 0, 1}-2g_{5, 0, 2}, \\
\medskip
g_{5, 1, 2}=\big(3r^2-9r+6\big)g_{3, 0, 1}-4g_{5, 0, 2}, \\
\medskip
g_{5, 2, 1}=\Big(\frac{nr}{2}-4r+4\Big)g_{4, 0, 1}+g_{5, 0, 2}, \\
\medskip
g_{5, 2, 2}=\frac{n(r-1)^4r}{2}+\big(-9r^2+24r-15\big)g_{3, 0, 1}+(-8r+12)g_{4, 0, 1}-5g_{5, 0, 1}
+6g_{5, 0, 2}, \\
\medskip
g_{5, 2, 3}=\Big(\frac{nr^2-9r^2-nr+21r-12}{2}\Big)g_{3, 0, 1}+2g_{5, 0, 2}, \\
\medskip
g_{5, 2, 4}=\big(3r^2-12r+12\big)g_{3, 0, 1}-2g_{5, 0, 2}, \\
\medskip
g_{5, 2, 5}=\Big(\frac{3r^2-15r+18}{2}\Big)g_{3, 0, 1}, \\
\medskip
g_{5, 3, 1}=\frac{n(r-2)(r-1)^3r}{2}+\big(-9r^2+33r-30\big)g_{3, 0, 1}+(-8r+16)g_{4, 0, 1}+4g_{5, 0, 2}, \\
\medskip
g_{5, 3, 2}=\frac{n(r-2)(r-1)^3r}{2}+\big(-12r^2+39r-30\big)g_{3, 0, 1}+(-4r+8)g_{4, 0, 1}+8g_{5, 0, 2}, \\
\medskip
g_{5, 3, 3}=\Big(\frac{3nr^2}{2}-12r^2-3nr+36r-24\Big)g_{3, 0, 1}+4g_{5, 0, 2}, \\
\medskip
g_{5, 4, 1}=\frac{n(r-1)^3r(nr-10r+8)}{4}+\Big(-3nr^2+39r^2+\frac{9nr}{2}-99r+60\Big)g_{3, 0, 1}+
\big(-2nr+28r-32\big)g_{4, 0, 1}+5g_{5, 0, 1}-14g_{5, 0, 2}, \\
\medskip
g_{5, 4, 2}=\frac{n(r-1)^3r(nr-12r+12)}{4}+\Big(\frac{-3nr^2+78r^2+3nr-204r+132}{2}\Big)g_{3, 0, 1}+\big(16r-24\big)g_{4, 0, 1}+5g_{5, 0, 1}-14g_{5, 0, 2}, \\
\medskip
g_{5, 4, 3}=\frac{n(r-2)^2(r-1)^2r}{8}+\big(-3r^2+12r-12\big)g_{3, 0, 1}+g_{5, 0, 2}, \\
\medskip
g_{5, 4, 4}=\frac{n(r-3)(r-2)(r-1)^2r}{6}+\big(-3r^2+15r-18\big)g_{3, 0, 1}, \\
\medskip
g_{5, 4, 5}=\Big(\frac{n^2r^2-16nr^2+72r^2+14nr-144r+72}{8}\Big)g_{3, 0, 1}-2g_{5, 0, 2}, \\
\medskip
g_{5, 5, 1}=\frac{n(r-2)(r-1)^2r(nr-10r+8)}{4}+\big(-3nr^2+42r^2+6nr-132r+96\big)g_{3, 0, 1}+(8r-16)g_{4, 0, 1}-12g_{5, 0, 2}, \\
\medskip
g_{5, 5, 2}=\frac{n(r-2)(r-1)^2r(nr-12r+12)}{12}+\big(9r^2-33r+30\big)g_{3, 0, 1}+(4r-8)g_{4, 0, 1}-2g_{5, 0, 2}, \\
\medskip
g_{5, 5, 3}=\frac{n(r-4)(r-3)(r-2)(r-1)r}{120}, \\
\medskip
g_{5, 6, 1}=\frac{n(r-3)(r-2)(r-1)r(nr-10r+8)}{48}+\Big(\frac{6r^2-30r+36}{4}\Big)g_{3, 0, 1}, \\
\medskip
g_{5, 6, 2}=\frac{n(r-1)^2r(n^2r^2-20nr^2+112r^2+14nr-176r+72)}{16}+\Big(\frac{-3n^2r^2+84nr^2-696r^2-90nr+1584r-888}{8}\Big)g_{3, 0, 1}\\
\medskip
~~~~~~~+\big(2nr-32r+36\big)g_{4, 0, 1}-5g_{5, 0, 1}+20g_{5, 0, 2}, \\
\medskip
g_{5, 6, 3}=\frac{n(r-1)^2r(n^2r^2-21nr^2+126r^2+16nr-216r+96)}{16}+\Big(\frac{18nr^2-252r^2-24nr+714r-444}{4}\Big)g_{3, 0, 1}+\big(nr-28r+36\big)g_{4, 0, 1}\\
\medskip
~~~~~~~-5g_{5, 0, 1}+17g_{5, 0, 2}, \\
\medskip
g_{5, 7, 1}=\frac{n(r-2)(r-1)r(n^2r^2-20nr^2+112r^2+14nr-176r+72)}{48}+\Big(\frac{3}{2}nr^2-21r^2-3nr+66r-48\Big)g_{3, 0, 1}+(-4r+8)g_{4, 0, 1}+4g_{5, 0, 2}, \\
\medskip
g_{5, 8, 1}=\frac{n(r-1)r(n^3r^3-30n^2r^3+328nr^3-1344r^3+18n^2r^2-444nr^2+3072r^2+160nr-2448r+672)}
{96}\\
\medskip
~~~~~~~+\Big(\frac{3n^2r^2-100nr^2+888r^2+106nr-1968r+1080}{8}\Big)g_{3, 0, 1}+\big(-2nr+36r-40\big)g_{4, 0, 1}+5g_{5, 0, 1}-20g_{5, 0, 2}, \\
\medskip
g_{5, 10, 1}=\frac{nr}{3840}\Big(n^4r^4-40n^3r^4+640n^2r^4-4960nr^4+16128r^4+20n^3r^3-720n^2r^3+9440nr^3-46080r^3+220n^2
r^2-6400nr^2\\
\medskip
\medskip
\medskip
~~~~~~~+51840r^2+1520nr-26880r+5376\Big)+\Big(\frac{-5n^2r^2+140nr^2-1080r^2-130nr+2160r-1080}{40}\Big)g_{3, 0, 1}+\Big(\frac{nr}{2}-8r+8\Big)g_{4, 0, 1}-g_{5, 0, 1}+4g_{5, 0, 2}. \\
\end{array}$
\endgroup
\begin{figure}[h]\label{sub}
\hspace{2.2cm} 
	\includegraphics[width=9cm,height=10.9cm]{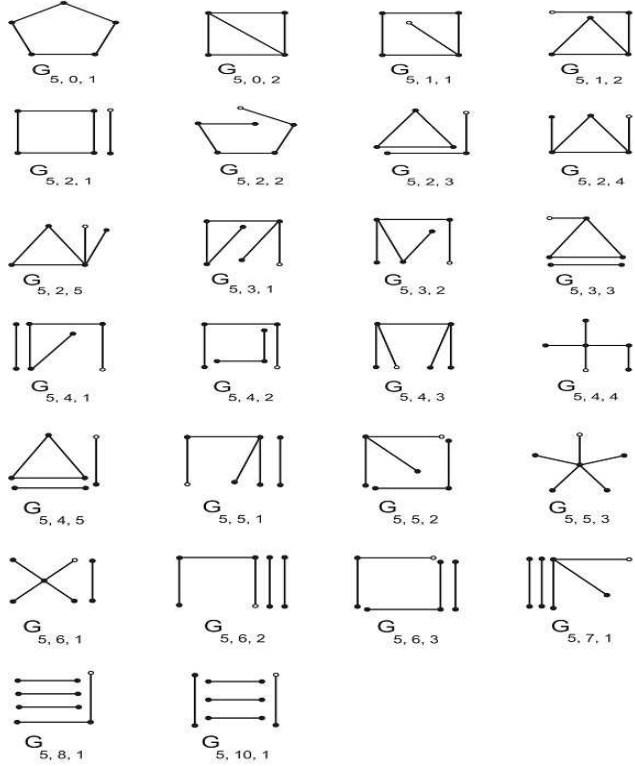}
 \caption{\label{figure 2}\small {The possible subgraphs with $5$ edges. The vertex $w$ is adjacent to open-circle vertex}.}
\end{figure}

\begin{figure}[!h]
	\centerline{\includegraphics[width=1cm]{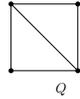}}
	\caption{\label{figure Q}\small {The subgraph $Q$. }}
\end{figure}

\medskip
The following theorem gives the matching polynomial of cubic graphs of order $10$ as a specific kind of regular graphs. 

Note that for every cubic graph $G$ of order $10$, $\alpha^{\prime}(G)=5$ and so it is enough to obtain 
$g_{m, 2m, 1}$ for every $m\leq 5$.

\begin{theorem}\label{pol}
Let $G$ be a cubic graph of order $10$. Then the matching polynomial of $G$ is 
\begin{flalign*}
\mu(G, x)&=x^{10}-15x^8+75x^6-145x^4+90x^2-18+(x^4-3x^2+3)g_{3, 0, 1} \\
&\qquad+(x^2+1)g_{4, 0, 1}+g_{5, 0, 1}-4g_{5, 0, 2}
\end{flalign*}
where $g_{3, 0, 1}$, $g_{4, 0, 1}$ and $g_{5, 0, 1}$ are the number of circuits on $3$, $4$ and $5$ vertices, respectively and $g_{5, 0, 2}$ is the number of subgraph $Q$, where $Q$ is the graph in Figure \ref{figure Q}. 
\end{theorem}

\proof
It is  sufficient to compute the coefficients of matching polynomial of $G$. We apply  the formulas for $g_{m, 2m, 1}$ where $0\leq m\leq 5$. We have $g_{1, 2, 1}=15$, $g_{2, 4, 1}=75$, $g_{3, 6, 1}=145-g_{3, 0, 1}$, $g_{4, 8, 1}=90-3g_{3, 0, 1}+g_{4, 0, 1}$ and $$g_{5, 10, 1}=18-3g_{3, 0, 1}-g_{4, 0, 1}-g_{5, 0, 1}+4g_{5, 0, 2}.$$ Now the result follows.
\qed

\medskip
\medskip
Behmaram proved that triangular-free $3$-regular graphs with ten vertices are matching unique and so cubic graphs $G_6, G_7, G_8, G_{10}, G_{16}$ and $G_{17}$ are matching unique (\cite{Beh}).
In the following theorem, we show that every cubic graph of order $10$ is matching unique.

\begin{theorem}
	Every cubic graph of order  ten  is matching unique.
\end{theorem}

\proof
Suppose that  $G_1, G_2, ..., G_{21}$ denote the cubic graphs of order $10$ as shown in Figure 1.
By replacing $g_{3, 0, 1}$, $g_{4, 0, 1}$, $g_{5, 0, 1}$ and $g_{5, 0, 2}$ for the cubic graphs with ten vertices in Theorem \ref{pol}, 
we have

$$\mu(G_1, x)=x^{10}-15x^8+75x^6-141x^4+80x^2-8,$$
$$\mu(G_2, x)=x^{10}-15x^8+75x^6-143x^4+86x^2-6,$$
$$\mu(G_3, x)=x^{10}-15x^8+75x^6-142x^4+84x^2-7,$$
$$\mu(G_4, x)=x^{10}-15x^8+75x^6-143x^4+90x^2-10,$$
$$\mu(G_5, x)=x^{10}-15x^8+75x^6-141x^4+80x^2-6,$$
$$\mu(G_6, x)=x^{10}-15x^8+75x^6-145x^4+95x^2-13,$$
$$\mu(G_7, x)=x^{10}-15x^8+75x^6-145x^4+93x^2-9,$$
$$\mu(G_8, x)=x^{10}-15x^8+75x^6-145x^4+92x^2-8,$$
$$\mu(G_9, x)=x^{10}-15x^8+75x^6-144x^4+89x^2-8,$$
$$\mu(G_{10}, x)=x^{10}-15x^8+75x^6-145x^4+95x^2-11,$$
$$\mu(G_{11}, x)=x^{10}-15x^8+75x^6-143x^4+87x^2-7,$$
$$\mu(G_{12}, x)=x^{10}-15x^8+75x^6-143x^4+84x^2-6,$$
$$\mu(G_{13}, x)=x^{10}-15x^8+75x^6-144x^4+91x^2-8,$$
$$\mu(G_{14}, x)=x^{10}-15x^8+75x^6-142x^4+81x^2-6,$$
$$\mu(G_{15}, x)=x^{10}-15x^8+75x^6-144x^4+90x^2-9,$$
$$\mu(G_{16}, x)=x^{10}-15x^8+75x^6-145x^4+96x^2-12,$$
$$\mu(G_{17}, x)=x^{10}-15x^8+75x^6-145x^4+90x^2-6,$$
$$\mu(G_{18}, x)=x^{10}-15x^8+75x^6-141x^4+84x^2-4,$$
$$\mu(G_{19}, x)=x^{10}-15x^8+75x^6-143x^4+85x^2-7,$$
$$\mu(G_{20}, x)=x^{10}-15x^8+75x^6-139x^4+78x^2-12,$$
$$\mu(G_{21}, x)=x^{10}-15x^8+75x^6-141x^4+90x^2-18.$$

\medskip
Assume that $H$  is a cubic graph with $10$ vertices. By the matching polynomials of $G_1, G_2, ..., G_{21}$ we have the result.
\qed




\begin{thebibliography}{1}
	
	
	\bibitem{turk} S. Alikhani and Y.H. Peng, {\it Domination polynomials of cubic graphs of order $10$}, Turk. J. Math.
	35 (2011) 355-366.
	
	\bibitem{jas} S. Alikhani and S. Jahari, {\it On the edge cover polynomial of certain graphs}, J. Alg. Sys. 2 (2) (2014) 97-108. 
	
	\bibitem{Ves} V. Andova, F. Kardo\v s and R. \v Skrekovski, {\it Sandwiching saturation number of fullerene graphs}, arXiv:1405.2197 (2014).

           \bibitem{Bee} R. A. Beezer, {\it The number of subgraphs of a regular graph}, Congressus Numerantium 100 (1994) 
89-96.

           \bibitem{Bze} R.A. Beezer, {\it Counting subgraphs in regular graphs}, Available at \texttt{http://buzzard.ups.edu/talks/beezer-2006-uwt-counting-subgraphs.pdf}.

           \bibitem{bee2} R.A. Beezer and E. J. Farrell, {\it The matching polynomial of a regular graph}, Discrete Math. 137 
(1995) 7-18.  

          \bibitem{Beh} A. Behmaram, {\it On the number of 4-matchings in graphs}, MATCH Commun. Math. Comput. Chem. 62 (2009) 381-388.

           \bibitem{Bie} T. Biedl, E. D. Demaine, C. A. Duncan, R. Fleischer and S. G. Kobourov, {\it Tight bounds on maximal and maximum matchings}, Discrete Math. 285.1 (2004) 7-15.


	 

	 \bibitem{Fad} J. Faudree, R. J. Faudree, R. J. Gould and M. S. Jacobson, {\it Saturation numbers for trees}, Electron. J. Combin. 16 (2009).
	 
	\bibitem{Gutman1}  I. Gutman, {\it The matching polynomial}, MATCH Commun. Math. Comput. Chem., 6 (1979), 	  75--91.
	  
	  \bibitem{Gutman2}  I. Gutman, {\it Uniqueness of the matching polynomial}, MATCH Commun. Math. Comput.	  Chem., 55 (2006), 351--358.
	  
	  \bibitem{Gutman3} I. Gutman, {\it Characteristic and matching polynomials of benzenoid hydrocarbons}, J. Chem.
	  Soc., Faraday Trans. 2., 79 (1983), 337--345.
	  
	  \bibitem{Gutman4}  I. Gutman and F. Harary, {\it Generalizations of the matching polynomial}, Utilitas Math., 24 	  (1983), 97--106.

	
        \bibitem{Kho} G.B. Khosrovshahi, Ch. Maysoori and B. Tayfeh-Rezaie, {\it A Note on $3$-Factorizations of 
$K_{10}$}, J. Combin. Designs 9.5 (2001) 379-383.

\bibitem{Match} R. Vesalian, F.Asgari, {\it Number of $5$-Matchings in Graphs}, MATCH Commun. Math. Comput. Chem. 69 (2013) 33--46. 

\end{thebibliography}
\end{document}